\documentstyle[amsfonts,12pt]{article}
\input amssym.def
\newcommand{\al}{{\alpha}}
\newcommand{\la}{{\lambda}}

\newcommand{\Si}{{\Sigma}}

\newcommand{\ga}{{\gamma}}
\newcommand{\vf}{{\varphi}}
\newcommand{\ka}{{\kappa}}
\newcommand{\G}{{\Gamma}}

\newcommand{\pa}{{\partial}}

\newcommand{\cP}{{\cal P}}

\newcommand{\cB}{{\cal B}}

\newcommand{\bbP}{{\Bbb P}}
\newcommand{\bbC}{{\Bbb C}}
\newcommand{\bbZ}{{\Bbb Z}}

\newcommand{\te}{\theta}

\newcommand{\de}{{\delta}}
\newcommand{\D}{{\Delta}}

\newcommand{\tL}{{\tilde{L}}}
\date{}

\newtheorem{lem}{Lemma}
\newtheorem{teo}{Theorem}

\begin{document}

\vspace{0.3in}
\begin{flushright}
ITEP-TH-44/02\\
math.AG/0209145\\
\end{flushright}
\vspace{10mm}

\begin{center}
{\Large\bf R-Matrix Structure of Hitchin System \\[0.3cm]
in Tyurin Parameterization}\\[0.3cm]
V.A. Dolgushev\footnote{E-mail: vald@mit.edu}\\[0.3cm]
{\it Department of Mathematics, MIT,} \\
{\it 77 Massachusetts Avenue,} \\
{\it Cambridge, MA, USA 02139-4307;} \\
{\it University Center, JINR, Dubna,} \\
{\it 141 980, Moscow Region, Russia;} \\
{\it Institute for Theoretical and Experimental Physics, } \\
{\it 117259, Moscow, Russia.}
\end{center}

\begin{abstract}
We present a classical $r$-matrix for the Hitchin
system without marked points on an arbitrary
non-degenerate algebraic curve of genus $g\ge2$ using
Tyurin parameterization of holomorphic vector
bundles.
\end{abstract}

\section{Introduction}
The study of the moduli space of holomorphic
vector bundles over an algebraic curve motivated
by the geometric Langlands conjecture \cite{BD}
is now one of the most fascinating topics of modern
algebraic geometry. Important tools for
the investigation are integrable systems
of Hitchin type \cite{Hit}, \cite{IM}, \cite{LOZ},
\cite{Mark}, \cite{Nikita} whose configuration spaces are defined as
connected components of the moduli space
of holomorphic vector bundles over
compact Riemann surfaces.

It is not an easy task to give a satisfactory description of a
Hitchin system since its definition is implicit and
at the first sight it is clear neither how to find
its Lax representation nor how to write down
the respective equations of motion. For the case of algebraic curves
of genus zero and one this question was solved
in papers \cite{ER1}, \cite{FFR}, \cite{Felder},
\cite{Nikita} and for Schottky curves of an arbitrary
higher genus a description of Hitchin systems
was proposed in \cite{Schottky}.

In the paper by A. Tyurin \cite{AT} a classification of holomorphic vector
bundles over algebraic curves of arbitrary genus is obtained and a convenient
parameterization of big cells of connected components of the moduli space of
the bundles is suggested. After its introduction in \cite{AT}
Tyurin parameterization is used in the works on integrable
differential systems \cite{ER},  \cite{A2}, \cite{IM}, \cite{A1}, \cite{A3}
 and two papers \cite{ER}, \cite{IM} are worthy of
mention, in which the Tyurin description is used to parameterize
Hitchin systems.

In \cite{ER} the parameterization of Hitchin systems
is obtained for the case of rank $2$ holomorphic vector
bundles of degree $2g$ over algebraic curves of genus $g\ge 1$
and in \cite{IM}
Tyurin description is used to parameterize
an arbitrary Hitchin system
and to construct infinite-dimensional field analogues\footnote{In this context
paper \cite{LOZ} is also worthy of mention, in which the case of
two-dimensional version of the elliptic Gaudin system is
considered in detail.} of the systems of Hitchin type.
The results of the papers \cite{ER} and
 \cite{IM} show that
Tyurin parameterization should enable us to achieve a
rough but explicit description of quantum Hitchin
systems, and the first step in this direction is
a construction of classical $r$-matrix structures
for Hitchin systems which will allow us
to quantize the systems in a quantum
group theoretic setting \cite{EV}, \cite{DQG1}, \cite{DQG}.

The concept of a classical $r$-matrix was originally introduced in
works of ``Leningrad school'' \cite{STS}, \cite{Sk79}
(see also book \cite{FT}) as a natural object that
encodes the Hamiltonian structure of the Lax equation,
provides the involution of integrals of motion \cite{BV90},
and gives a natural framework for quantizing integrable systems.

In this paper we present a classical $r$-matrix for
the Hitchin system without marked points on an arbitrary
non-degenerate algebraic curve of genus $g\ge2$ using
Tyurin parameterization of the moduli space
of rank $n$ holomorphic vector bundles
of degree $ng$.

Following Tyurin
\cite{AT} a generic holomorphic vector bundle
$\cB$ of this type over a non-degenerate curve $\Si$
has an $n$-dimensional space $H^0(\Si,\cB)$ of
holomorphic sections and for a
generic point $P$ of the curve $\Si$ the sections generate
a basis in the respective fiber $\cB_P$. However, the evaluations
of these sections on $ng$ points
$\ga_a\in \Si$, $a=1,\ldots, ng$ are
linearly dependent and for distinct points $\ga_a$
they determine subspaces $V_a\in \cB_{\ga_a}$ of
codimension $1$ or just one-dimensional
linear subspaces $l_a$ in the dual
space $H^0(\Si,\cB)^*$. The collection
of lines $l_a\in H^0(\Si,\cB)^*$ can be identified
with nonzero vectors $\al_a\in \bbC^n$, which are
defined up to the scalar multiples

$$
\al_a \mapsto \la_a \al_a, \qquad \la_a\in \bbC,~\la_a\neq 0
$$
and up to the following transformations of the group
$SL_n(\bbC)$\footnote{We assume here the summation
over repeated indices.}

\begin{equation}
\al^i_a\mapsto \tilde\al^i_a=\al^j_a (G^{-1})^j_i, \qquad
\det||G^i_j||=1,~i,j=1,\ldots n,
\label{diagact}
\end{equation}
generated by changes of basis in $H^0(\Si,\cB)$.

Thus, we arrive at the Tyurin map from an open dense
set of the moduli space
of rank $n$ holomorphic vector bundles
of degree $ng$ over the curve $\Si$ to
the following quotient

\begin{equation}
[\Si\times \bbP(\bbC^n)]^{(ng)}/SL_n(\bbC),
\label{cell}
\end{equation}
where the notation $^{(ng)}$ stands for the
symmetric direct product.

The set of points $\ga_a\in\Si$ and vectors
$\al_a\in \bbC^n$ are referred to as Tyurin parameters,
and the main statement of the paper \cite{AT}
we are going to use is that an open dense set of the
moduli space of rank $n$ holomorphic
vector bundles of degree $ng$ over the curve
$\Si$ is parameterized by points of
the quotient (\ref{cell}).

Note that in our considerations we are not going
to bother about the singularities of the quotient of
the space $[\Si\times \bbP(\bbC^n)]^{ng}$
with respect to the action of the
symmetric group $S_{ng}$
and in what follows
we omit factorization of our
phase space with respect to
permutations.

In order to parameterize
the phase space of the Hitchin system without marked
points one has to supplement Tyurin parameters
$(\ga_a,\al_a)$ with points $\ka_a\in T^*_{\ga_a}\Si$
and vectors $\beta_a\in \bbC^n$, which are subject to the
following conditions

\begin{equation}
\sum^n_{i=1}\beta^i_a\al^i_a=0,
\label{chastInt}
\end{equation}

\begin{equation}
T_{ij}=\sum^{ng}_{a=1}\beta^i_a\al^j_a= 0.
\label{vychet}
\end{equation}

Equation (\ref{chastInt}) means that the $\beta^i_a$
may be regarded as coordinates in the
cotangent space $T^*_{\al_a}\bbP^{n-1}$ and
equations (\ref{vychet}) are just the first class
constraint conditions corresponding to the symplectic
action of the group $SL_n(\bbC)$ on the parameters
$\al^i_a$ and $\beta^i_a$\,, that is

\begin{equation}
\al^i_a\mapsto \al^j_a(G^{-1})^j_i, \qquad
\beta^i_a \mapsto G^i_j \beta^j_a, \qquad
||G^i_j||\in\, SL_n(\bbC).
\label{action}
\end{equation}

In other words, the phase space of the
Hitchin system in question can obtained via
symplectic reduction in the space

\begin{equation}
\cP=T^*[\Si\times \bbP(\bbC^n)]^{ng}
\label{Shirokoe}
\end{equation}
on the surface of the first class
constraints\footnote{A similar trick is used in \cite{ER1} for
description of Hitchin systems associated with
marked rational and elliptic curves.}
(\ref{vychet}).

The main statement of the paper (see Theorem 1) is
that the Krichever Lax matrix of the
Hitchin system being extended to the symplectic manifold
(\ref{Shirokoe}) admits a simple $r$-matrix
structure, which is defined by a matrix-valued meromorphic
section of the bundle
$\Si\times T^*\Si$ over the direct product
of curves $\Si\times \Si$. We argue that
using the $r$-matrix structure
one can easily derive the classical $r$-matrix for the initial
Lax matrix of the Hitchin system either with a help of
a gauge invariant extension of
the Krichever Lax matrix to the manifold
(\ref{Shirokoe}) or with the help of
on-shell Dirac brackets between the entries of
the initial extension of the Krichever Lax matrix.

Note however that
the $r$-matrix structure of the extended system
is much simpler than the resulting $r$-matrix
of the Hitchin system and this remarkable
simplification turns out to be possible due
to the fact\footnote{I am indebted to A.M. Levin
for the technical trick concerning
the extension of the Krichever Lax matrix.} that the Krichever Lax matrix
of the Hitchin system, being a meromorphic
differential on the curve $\Si$\,, can be extended to the
symplectic manifold (\ref{Shirokoe}) in such a way
that the extension is also a meromorphic matrix-valued
differential on $\Si$.

The organization of the paper is as follows.
In the second section we present the extension of
the Krichever Lax matrix for the Hitchin system without marked
points on a non-degenerate algebraic curve of
genus $g\ge 2$ and propose that the extended system
admits an $r$-matrix structure, which is defined
as a meromorphic matrix-valued function
on one copy of the curve and a meromorphic
$1$-form on another copy of the same curve.
Then, postponing the proof of this proposition
to the next section, we show how to derive
the classical $r$-matrix for the genuine
Krichever Lax matrix of the Hitchin system using
the above $r$-matrix structure.

Before presenting the proof in section 3 we show
that a matrix-valued differential that enters into the
definition of the above $r$-matrix structure does exist.
We also give the properties of the differential
as a function in the first variable and identify
derivatives of the extended Krichever Lax matrix with respect
to phase space variables as meromorphic
differentials on the algebraic curve.

In the concluding section we mention dynamical properties
of the presented $r$-matrices, discuss
a possibility to derive the $r$-matrices using an
infinite-dimensional Hamiltonian reduction, and raise some
other questions.

In the appendix at the end of the paper we
present the Krichever lemma, which is used
throughout the paper as a tool, that enables us
to identify meromorphic vector-valued differentials
by their singular parts and certain linear
equations for their regular parts.
Although the statement is analogous to
Lemma $2.2$ in \cite{IM} we present its proof
in the appendix since in some respect
the statement generalizes the
lemma and the presented proof differs from
the one given in \cite{IM}.

In this paper we use
standard notations for Poisson brackets
between entries of a Lax matrix.
For example, if
$L(z)$ is a matrix-valued function and

\begin{equation}
r(z,w)=\sum_{i,j,k,l}r_{ijkl}(z,w)e_{ij}\otimes e_{kl},
\label{RMnotta}
\end{equation}
where

$$(e_{ij})_{kl}=\de_{ik}\de_{jl}$$
are the elements of the standard basis in $gl_n(\bbC)$
then the expression
$$\{L_1(z),L_2(w)\}=[r(z,w), L_1(z)]
-[r_{21}(w,z), L_2(w)]$$
means that Poisson brackets between
the entries $L_{ij}(z)$ and $L_{kl}(w)$
take the following form

$$
\{L_{ij}(z),L_{kl}(w)\}=\sum_{m=1}^n (r_{imkl}(z,w) L_{mj}(z)-
L_{im}(z)r_{mjkl}(z,w))-
$$
$$
-\sum_{m=1}^n (r_{kmij}(w,z) L_{ml}(w)-
L_{km}(w)r_{mlij}(w,z)).
$$

Throughout the paper we assume that
$\Si$ is a non-degenerate algebraic curve
of genus $g\ge 2$.

\section{$R$-matrix structure for the Hitchin system
without marked points.}

We start with the following particular case
of Lemma $2.2$ in \cite{IM}

\begin{lem} For a generic set of pairs
$(\ga_a, k_a),$ $\ga_a\in \Si,$ $k_a\in T^*_{\ga_a}\Si,$
$a=1,\ldots,ng$ and complex parameters $\al^i_a$
and $\beta^i_a$ $i=1,\ldots n$
such that

\begin{equation}
\sum^n_{i=1}\beta^i_a\al^i_a=0
\label{chast}
\end{equation}
there exists a unique matrix-valued
meromorphic differential $L_{ij}=L_{ij}(z)dz$
of the third kind satisfying
the following properties

\begin{enumerate}
\item The differential $L_{ij}$ has
poles only at the points $\ga_a$ and at some
fixed point $P\in \Si$.

\item On a neighborhood of the point $\ga_a$ the
differential $L_{ij}(z)dz$ behaves like

\begin{equation}
L_{ij}(z)=\frac{\beta^i_a\al^j_a}{z-z(\ga_a)}+
L_{ij}^{a,0}+ L_{ij}^{a,1} (z-z(\ga_a))+\ldots.
\label{LAX}
\end{equation}

\item $\al_a$ is a left eigenvector
for the matrix $||L_{ij}^{a,0}||$ with the
eigenvalue $\ka_a$

\begin{equation}
\sum_{i=1}^n \al_a^i L_{ij}^{a,0} =\ka_a \al^j_a.
\label{eigen}
\end{equation}
\end{enumerate}
\end{lem}

The differential $L_{ij}(z)dz$ is obviously
invariant under the transformations

\begin{equation}
\al_a\mapsto \la_a \al_a, \qquad
\beta_a\mapsto \la^{-1}_a \beta_a, \qquad
\la_a\in \bbC, \quad \la_a\neq 0
\label{equiv}
\end{equation}
and, hence, it may be regarded as a function
with values in meromorphic differentials
on an open dense set of the space (\ref{Shirokoe})
so that the components of the
vector $\al_a$ are identified
with homogeneous coordinates
in $\bbP(\bbC^n)$ and the components of
the vector $\beta_a$, being subject to
the conditions (\ref{chast}) define a
point in the respective cotangent space
$T^*_{\al_a}\bbP(\bbC^n)$.

The differential $L_{ij}(z)dz$ and its natural generalizations
were originally found in the paper \cite{IM}
by Krichever as solutions of the momentum map
equations for Hitchin systems. Although the differential
$L_{ij}(z)dz$ is not a Krichever Lax matrix of the
Hitchin system without marked points since
equations (\ref{vychet}) are not imposed,
$L_{ij}(z)dz$ may be regarded as
an extension of the above Lax matrix to the
symplectic manifold $\cP$. In what follows, we
refer to $L$ as a Krichever Lax differential.

Notice that in view of Lemma $2.1$ of \cite{IM}\,,
the differential $L_{ij}(z)dz$ can
be identified with a meromorphic section with
a single pole at the point $P$ of the bundle
$End(\cB)\otimes {\cal K}$ where $\cB$ is the holomorphic bundle over
$\Si$ corresponding to the Tyurin parameters $\ga_a$ and $\al^i_a$ and
${\cal K}$ is a canonical bundle of the curve $\Si$.

Soon we will show that the Krichever Lax differential being considered as a function
on the symplectic manifold $\cP$ admits an $r$-matrix structure, but  now we
present an important ingredient which enters into the definition of the $r$-matrix
structure in question.

\begin{lem}
\label{LEMMA}
For a generic set of Tyurin parameters $\al_a$ and
$\ga_a$ there exists a unique matrix-valued differential
$r_{jk}(z,w)dw$ such that

\begin{enumerate}
\item  $r_{jk}(z,w)dw$ is a meromorphic
function in $z$ and a meromorphic $1$-form in $w$,

\item $r_{jk}(z,w)dw$ is holomorphic in $w$
everywhere on $\Si$ except the points $w=w(P)$ and $w=z$, where
it has simple poles with residues
$\de_{jk}$ and $-\de_{jk}$\,, respectively,

\item $\al_{a}$ are null vectors for
the matrices $r_{jk}(z,\ga_a)$

\begin{equation}
\sum_{k=1}^n r_{jk}(z,\ga_a)\al^k_a=0.
\label{nullvv}
\end{equation}
\end{enumerate}
\end{lem}

The existence of the meromorphic differential $r_{jk}(z,w)dw$
which is also a meromorphic function in $z$ satisfying the
above conditions is proved in subsection $3.1$ where
a stronger statement (see Lemma \ref{dalshe})
concerning the properties of the differential $r_{jk}(z,w)dw$ as
a function in $z$ is also formulated.
To this end, the uniqueness of the differential
$r_{jk}(z,w)dw$ follows directly from the Krichever
lemma.

We now present the main statement of the paper.

\begin{teo}
For an arbitrary non-degenerate algebraic curve $\Si$
of genus $g\ge 2$ the canonical
Poisson brackets of the space (\ref{Shirokoe})
between the entries of
the Krichever Lax differential (\ref{LAX}) obey the Yang-Baxter
relation

\begin{equation}
\{L_1(z),L_2(w)\}dz\otimes dw=[r(z,w), L_1(z)]dz\otimes dw
-[r_{21}(w,z), L_2(w)]dz\otimes dw,
\label{LrL}
\end{equation}
where the differential $r(z,w)dw$
is given by the formula

\begin{equation}
r(z,w)dw=\sum_{i,j,k}r_{jk}(z,w)e_{ij}\otimes e_{ki}dw,
\label{simp}
\end{equation}
and $r_{jk}(z,w)dw$ is the meromorphic
$1$-form defined in Lemma \ref{LEMMA}.
\end{teo}

We will refer to the differential (\ref{simp}) as
an $r$-matrix differential.

In the following section we present
an algebraic-geometric proof of the theorem.
First we explain how to achieve the $r$-matrix
for the Hitchin system we consider
using the differential (\ref{simp}).

As we have mentioned in the introduction,
the phase space of the Hitchin system without marked
points can be identified with an open dense set of the
quotient of the constraint surface (\ref{vychet})
in the space (\ref{Shirokoe}) with respect to the
symplectic action (\ref{action})
of the group $SL_n(\bbC)$.

In other words, if one chooses some gauge fixing
conditions

\begin{equation}
\chi^{ij}(\al^k_a)=0
\label{chi}
\end{equation}
for the transformations
(\ref{action}) then the phase space of the Hitchin system
can be roughly identified with an intersection
of the surfaces (\ref{vychet}) and (\ref{chi}) in the
space (\ref{Shirokoe}) and the respective Krichever Lax matrix
is defined as the differential
(\ref{LAX}), restricted to the intersection

\begin{equation}
l_{ij}(z)dz=L_{ij}(z)dz|_{T_{kl}=\chi^{kl}=0}.
\label{laxic}
\end{equation}

Obviously, the Lax matrix (\ref{laxic}) is a meromorphic differential on the
curve $\Si$ with the same properties (\ref{LAX}), (\ref{eigen}) as the Krichever Lax
differential except that the point $P$ is now regular for the differential
(\ref{laxic}). In view of Lemma $2.1$ of \cite{IM}, this means that the
differential (\ref{laxic}) can be identified with a holomorphic section of
the bundle $End(\cB)\otimes{\cal K}$.

The gauge transformations (\ref{action}) of the Krichever Lax differential (\ref{LAX})
have a form of adjoint action

\begin{equation}
L(z)\mapsto G L(z) G^{-1},
\qquad
G\in\, SL_n(\bbC),
\label{actionL}
\end{equation}
and hence, whatever gauge fixing conditions
are chosen, the desired $r$-matrix of the Hitchin system
can be derived from the $r$-matrix differential (\ref{simp})
either with the help of a gauge invariant extension of
the Krichever Lax matrix (\ref{laxic}) or with the help of
on-shell Dirac brackets between the entries of
the differential (\ref{LAX}).

Recall that the gauge invariant extension of
Lax matrices was originally used to calculate classical
$r$-matrices for integrable systems in the works
\cite{ArMed} and \cite{Babel}. In a more general situation
Dirac brackets and gauge invariant extension
of Lax matrices are used for analogous calculations
in the paper \cite{BDOZ}. At last, in the paper
\cite{Vengry} Dirac bracket technique is used in a
specific framework to obtain new examples of
Etingof-Varchenko dynamical $r$-matrices \cite{EVclass}.

To derive the classical $r$-matrix for the Hitchin system we use the gauge
invariant extension of the Krichever Lax matrix (\ref{laxic})
to the space (\ref{Shirokoe}).
For example,
if some $n\times n$-minor $||\al^i_{a_j}||$ of the matrix $||\al^i_a||$
where $1\le a_1<a_2<\ldots <a_n\le ng$ is non-degenerate we can choose gauge
fixing conditions in the form \cite{IM}

\begin{equation}
\al^i_{a_j}=0,~ i\neq j,
\qquad
\al^1_b=\al^2_b=\ldots =\al^n_b,
\label{gaugefix}
\end{equation}
where $b$ does not coincide with any of the
indices $a_1,\,a_2,\,\ldots \, a_n$.

On an open region of the space (\ref{Shirokoe}) one can define the
$SL_n$-valued function $G(\al_a)$ such that
if the vectors $\al_a$  do not satisfy the gauge fixing conditions
(\ref{gaugefix}) then the
transformed vectors $\tilde\al_a$
$$
\tilde\al^i_a=\al^j_a (G^{-1}(\al_c))^j_i
$$
do so. Otherwise $G(\al_a)$
is just the identity matrix.

Then the matrix-valued differential

\begin{equation}
l^G(z)=G(\al_a)L(z)G^{-1}(\al_a)
\label{gaugecont}
\end{equation}
turns out to be a desired gauge invariant extension
of the Krichever Lax matrix (\ref{laxic})
to the space $\cP$ and the
$r$-matrix in question
takes the form

\begin{equation}
r^{H}(z,w)dw=
(r(z,w)dw+\{G_1(\al_a),L_2(w)\}dw)|_{on~shell},
\label{viacont}
\end{equation}
where the notation $|_{on~shell}$ means that
the expression in the parenthesis is considered
on the surface of the constraints
(\ref{vychet}) and (\ref{gaugefix}).

~\\
{\bf Example.} Although Hitchin systems without marked points are
non-trivial only for algebraic curves of genus $g\ge 2$
the Krichever Lax differential (\ref{LAX}) and its $r$-matrix
structure (\ref{simp}) exist on an elliptic curve
as well.

To show this, we realize an elliptic curve $\G$ as a quotient
$\G=\bbC/\{1,\tau\}$, $Im\,\tau>0$ and denote the parameters
$\ga_a$ and $k_a$ by $q_a$ and $p_a$, respectively, where $a$ now runs
from $1$ to $n$. Then, the Krichever Lax differential (\ref{LAX}) and the
$r$-matrix differential (\ref{simp}) can be written in terms
of the standard $\te$-function as follows

\begin{equation}
\begin{array}{c}
\displaystyle
 L_{ij}(z)=\sum_{k,l=1}^n\pi^k_i\tL_{kl}(z)\al^j_l, \qquad \tL_{ii}=p_i, \\[0.5cm]
\displaystyle
\tL_{ij}(z)=\sum_{k=1}^n \al^k_i \beta^k_j~
\frac{\te(z-q_i)\te(z+q_i-q_j)\te(q_j)\te'(0)}{\te(z)
\te(z-q_j)\te(q_j-q_i)\te(q_i)},\quad i\neq j,
\end{array}
\label{tL}
\end{equation}

$$r(z,w)=\sum_{i,j=1}^n (E(z-w)+E(w)) e_{ij}\otimes e_{ji}-$$
\begin{equation}
-\sum_{i,j,k,a=1}^n \pi_k^a\al^j_a (E(z-q_a) + E(q_a)) e_{ij}\otimes e_{ki},
\label{onA}
\end{equation}
where $||\pi_i^j||$ is the inverse matrix to $||\al_k^l||$
$$
\sum_{k=1}^n \al_i^k\pi_k^j=\de_i^j,
$$
$$
\te(z)=\sum_{m\in\bbZ}
exp\,(\pi i \tau (m+1/2)^2 + 2\pi i (m+1/2)(z+1/2)),
$$
and
$$E(z)=\frac{\te'(z)}{\te(z)}.$$

To explain the relation of the Lax matrix (\ref{tL}) to
Lax representation of known integrable systems we
have to enlarge the phase space parameterized by
coordinates $q_a,\, p_a,\, \al_a^i,$ and $\beta^i_a$ with
some coadjoint orbit ${\cal O}$ of the group $SL_N$.
Symplectic reduction of this space to the
first class constraint surface
\begin{equation}
\sum_{a=1}^{ng}\beta^i_a\al^j_a+\eta_{ij}=0
\label{eeeee}
\end{equation}
leads us to the phase space and the
Lax matrix\footnote{We note that the Lax matrix
of the elliptic spin Calogero-Moser system was originally presented as
a meromorphic function on the elliptic curve in the
paper \cite{IM}} of the elliptic
spin Calogero-Moser system \cite{Nikita}, \cite{Olshanet}.
Here $\eta_{ij}$ denote conventional coordinates
on the coadjoint orbit ${\cal O}$\,.
If we now restrict ${\cal O}$ to be the
maximal coadjoint orbit we just get the
particular case of one marked point of the integrable system
considered in \cite{ER1}. The latter system is
now generally regarded as an elliptic Gaudin
system \cite{Nikita}, \cite{Dima}.

\section{The proof of the Yang-Baxter relation}
The proof of Theorem 1 is based on the observation
that both sides of equation (\ref{LrL}) satisfy the
same properties, which, in turn, uniquely define
them as meromorphic forms on the direct product of curves
$\Si\times \Si$. Namely, it turns out that both
sides of equation (\ref{LrL}) have
coincident singular parts while
their regular parts at the points $\ga_a$
obey the same linear inhomogeneous equations,
which  uniquely define the remaining
arbitrariness in the holomorphic parts due
to the Krichever lemma.

To calculate Poisson brackets between the
entries of the Krichever Lax differential we choose the local chart of
the space (\ref{Shirokoe}) where
\begin{equation}
\al^1_a=1,\qquad
\beta^1_a=-\sum_{\mu=2}^n\beta^{\mu}_a\al_a^{\mu}, \qquad
\forall\,a=1,\ldots, ng.
\label{chart}
\end{equation}

Note that although a choice of another
local affine chart affects intermediate calculations
the Poisson bracket
\begin{equation}
\{L_{ij}(z),L_{kl}(w)\}dz\otimes dw
\label{PBLL}
\end{equation}
is, in fact, ``a function'' on the space (\ref{Shirokoe}), and
therefore the properties of the expression (\ref{PBLL})
as a form on product of curves $\Si\times \Si$ do not
depend on the choice of local coordinates on $\cP$.

Throughout this section we also assume that some local
coordinates are chosen on neighborhoods of the points $\ga_a$ on
the curve $\Si$ and for simplicity we denote the coordinate
$z(\ga_a)$ by the same letter $\ga_a$.

\subsection{Properties of the $r$-matrix differential as a function
of the first argument.}

We start this subsection with the following
\begin{lem}
\label{dalshe}
The differential $r_{ij}(z,w)dw$
defined in Lemma \ref{LEMMA} exists
and is holomorphic in $z$ everywhere on $\Si$
except the points $\ga_a$, where it has
simple poles. The differential $r_{ij}(z,w)dw$
is also vanishing at the point $z=z(P)$

\begin{equation}
\label{vP0}
r_{ij}(z(P),w)dw=0.
\end{equation}

\end{lem}

~\\
{\bf Proof.} First, using the Krichever lemma we introduce
auxiliary holomorphic vector-valued differentials $u_{ai}(z)dz$,
which are uniquely defined by the following properties

\begin{equation}
\sum_{i=1}^n u_{ai}(\ga_b)\al_b^i=\de_{ab}.
\label{uia}
\end{equation}

Using standard arguments based
on the Kodaira-Nakano vanishing theorem and
GAGA principles one can easily show that
for an arbitrary point $Q\in\Si$ there
exists a matrix-valued differential $\D_{ij}(z,w)dw$, which is
holomorphic in $z$ on some open neighborhood $U_Q$ of the
point $Q$ and holomorphic in $w$ everywhere on $\Si$ except
the points $w=w(P)$ and $w=z$, where the differential
has simple poles with
residues $\de_{ij}$ and $-\de_{ij}$, respectively.

It is easy to see that the following
matrix-valued differential

\begin{equation}
\label{rDm}
r^{U_Q}_{ij}(z,w)dw=\D_{ij}(z,w)dw- \sum_{a,k}\D_{ik}(z,\ga_a)\al^k_a u_{aj}(w)dw
\end{equation}
is meromorphic in $z$ on the neighborhood $U_Q$ and
satisfies conditions $2$ and $3$ of Lemma
\ref{LEMMA}.

Since conditions $2$ and $3$ of Lemma
\ref{LEMMA} uniquely determine $r^{U_Q}_{ij}(z,w)dw$ as
a $1$-form in $w$ we can define the desired differential $r_{ij}(z,w)dw$
by its restrictions (\ref{rDm}) to the sets $U_Q$.

Equation (\ref{rDm}) also implies that the resulting
differential $r_{ij}(z,w)dw$ is holomorphic in $z$ everywhere on $\Si$
except the points $\ga_a$ and on the neighborhoods of the points
the differential behaves like

\begin{equation}
\label{kakprim}
r_{ij}(z,w)dw=-\frac{\al^i_a u_{aj}(w)dw}{z-\ga_a}+{\rm regular~terms}.
\end{equation}

Note also that as the differential $r_{ij}(z(P),w)dw$ is
holomorphic in $w$ everywhere on $\Si$ equations (\ref{nullvv})
imply that the differential is in fact vanishing due to the Krichever
lemma.

Thus, the statement is proved. $\Box$

In order to prove the Yang-Baxter relation we have to
identify the next two coefficients of the Laurent expansion
of the differential $r_{ij}(z,w)dw$ in the first variable $z$
around a point $\ga_a$. In the following lemma
we identify these coefficients as meromorphic
differentials on the curve $\Si$.

\begin{lem}
\label{poperv}
The expansion coefficients $r^{a,0}_{ij}(w)dw$ and $r^{a,1}_{ij}(w)dw$
of the Laurent series

\begin{equation}
r_{ij}(z,w)dw=-\frac{\al^i_a u_{aj}(w)dw}{z-\ga_a}+
r^{a,0}_{ij}(w)dw+(z-\ga_a)r^{a,1}_{ij}(w)dw+ o(z-\ga_a)
\label{Laurent}
\end{equation}
of the differential $r_{ij}(z,w)dw$
on a neighborhood $U_{\ga_a}$ of a point $\ga_a$ are uniquely
defined by the following
properties\footnote{Note that the uniqueness of the
differentials $r^{a,0}_{ij}(w)dw$ and
$r^{a,1}_{ij}(w)dw$  satisfying the presented properties
follows from the Krichever lemma.}:

\begin{enumerate}

\item The $1$-form $r^{a,0}_{ij}(w)dw$ is holomorphic
everywhere on $\Si$ except the points
$P$ and $\ga_a$, where it has simple poles
with residues $\de_{ij}$ and $-\de_{ij}$,
respectively.

\item For $b\neq a$\,, $\al_b$ is a null vector for
the matrix $||r^{a,0}_{ij}(\ga_b)||$
$$
\sum_{j=1}^n r^{a,0}_{ij}(\ga_b)\al^j_b =0, \qquad b\neq a,
$$
and $\al_a$ is a null vector for
the regular part  the matrix
$||r^{a,0}_{ij}(w)||$ at the point $w=\ga_a$
$$
\sum_{j=1}^n r^{a,0}_{ij}(w)\al^j_a|_{{\rm regular~part~at}~w=\ga_a} =0.
$$

\item The $1$-form $r^{a,1}_{ij}(w)dw$ has a single
pole at the point $\ga_a$ and on a neighborhood of
the point it behaves like

\begin{equation}
\label{ra1}
r^{a,1}_{ij}(w)=-\frac{\de_{ij}}{(w-\ga_a)^2}+ {\rm regular~terms}.
\end{equation}

\item  For $b\neq a$\,,
$\al_b$ is a null vector for the matrix
$||r^{a,1}_{ij}(\ga_b)||$\,:
$$
\sum_{j=1}^n r^{a,1}_{ij}(\ga_b)\al^j_b =0,  \qquad b\neq a,
$$
and $\al_a$ is a null vector for
the regular part of the matrix
$||r^{a,1}_{ij}(w)||$ at the point $w=\ga_a$\,:
$$
\sum_{j=1}^n r^{a,1}_{ij}(w)
\al^j_a|_{{\rm regular~part~at}~w=\ga_a} =0.
$$

\end{enumerate}

\end{lem}

~\\
{\bf Proof.} Applying the properties of the $1$-form $r_{ij}(z,w)dw$
(see Lemma \ref{LEMMA}) to the Laurent
expansion (\ref{Laurent}) we get that
outside the neighborhood $U_{\ga_a}$ the differential
$r^{a,0}_{ij}(w)dw$ has only a simple pole
at the point $w=w(P)$ with the residue $\de_{ij}$,
the differential $r^{a,1}_{ij}(w)dw$ is holomorphic in
the region $\Si\setminus U_{\ga_a}$, and for $b\neq a$ $\al_b$
is a right null vector for the
matrices $||r^{a,0}_{ij}(\ga_b)||$ and
$||r^{a,1}_{ij}(\ga_b)||$
$$
\sum_{j=1}^n r^{a,0}_{ij}(\ga_b)\al^j_b=0, \qquad
\sum_{j=1}^n r^{a,1}_{ij}(\ga_b)\al^j_b=0, \qquad b\neq a.
$$

The expansion (\ref{Laurent}) cannot be used for the case
when $w$ is on the neighborhood $U_{\ga_a}$ because
$r_{ij}(z,w)$ is irregular at the point $z=w$.

In order to cure the problem we consider the function

$$
\vf_{ij}(z)=r_{ij}(z,w)+\frac{\de_{ij}}{w-z},
$$
which is already holomorphic at the point $z=w$, and
therefore the Laurent expansion

\begin{equation}
\label{LLaurent}
\vf_{ij}(z)=-\frac{\al^i_a u_{aj}(w)dw}{z-\ga_a}+
(r^{a,0}_{ij}(w))+\frac{\de_{ij}}{w-\ga_a})+
\end{equation}
$$
+(r^{a,1}_{ij}(w)+\frac{\de_{ij}}{(w-\ga_a)^2})
(z-\ga_a)+ o(z-\ga_a),
$$
of the function
is convergent on the neighborhood $U_{\ga_a}$ even
in the case when the point $w$ is on
the neighborhood.

Hence, we can apply the remaining properties of the
differential $r_{ij}(z,w)dw$ to expansion (\ref{LLaurent})
and finally get that on the neighborhood $U_{\ga_a}$ the differentials
$r^{a,0}_{ij}(w)dw$ and $r^{a,1}_{ij}(w)dw$ behave like
$$
r^{a,0}_{ij}(w)dw=-\frac{\de_{ij}}{w-\ga_a}+{\rm regular~terms},
$$
$$
r^{a,1}_{ij}(w)dw=-\frac{\de_{ij}}{(w-\ga_a)^2}+{\rm regular~terms},
$$
and $\al_a$ is a right null vector for the regular parts
of the matrices $||r^{a,0}_{ij}(w)||$ and
$||r^{a,1}_{ij}(w)||$ at the point $w=\ga_a$.

Thus, the lemma is proved. $\Box$

\subsection{Derivatives of the Krichever Lax differential.}
In this subsection we
present the properties
of derivatives of the differential (\ref{LAX})
with respect to the variables $\ga_a$ and $\ka_a$ and
with respect to the canonical coordinates $\al^{\mu}_a$
and $\beta^{\mu}_a$ $\mu=2,\,\ldots,n$ in the local
chart (\ref{chart}) on the space $\cP$. As it will be seen
the properties uniquely define the derivatives of $L$
as meromorphic differentials
on the curve $\Si$.

First, we note that the differential $\pa_{k_a}L_{ij}(z)dz$
can be written in the following form

\begin{equation}
\label{pakL}
\pa_{k_a}L_{ij}(z)dz=\al_a^j u_{ai}(z) dz,
\end{equation}
where $u_{ai}(z)dz$ are holomorphic differentials defined
by equations (\ref{uia}).

Second, the differential $\pa_{\beta^{\mu}_a}L_{ij}(z)dz$ has
at most simple poles at the points $P$ and $\ga_a$ and
the residue of $\pa_{\beta^{\mu}_a}L_{ij}(z)dz$ at
the point $\ga_a$ equals

\begin{equation}
\label{pabeL}
Res_{z=\ga_a}\pa_{\beta^{\mu}_a}L_{ij}(z)dz=
\de_{i\mu}\al^j_a -\de_{i1}\al_a^{\mu}\al_a^j.
\end{equation}

For $b\neq a$\,, $\al_b$ is a left null vector for the matrix
$||\pa_{\beta^{\mu}_a}L_{ij}(\ga_b)||$
$$
\sum_{i=1}^n\al^i_b
\pa_{\beta^{\mu}_a}L_{ij}(\ga_b)=0,  \qquad b\neq a,
$$
and $\al_a$ is a left null vector for the
regular part of the matrix
$||\pa_{\beta^{\mu}_a}L_{ij}(z)||$ at point $\ga_a$\,:
$$
\sum_{i=1}^n
\al^i_a \pa_{\beta^{\mu}_a}L_{ij}(z)|_{{\rm regular~part~at}~z=\ga_a}=0.
$$

Third, the differential $\pa_{\al^{\mu}_a}L_{ij}(z)dz$ also
has at most simple poles at the points $P$ and $\ga_a$ and
the residue of $\pa_{\al^{\mu}_a}L_{ij}(z)dz$ at
the point $\ga_a$ equals

\begin{equation}
\label{paalL}
Res_{z=\ga_a}\pa_{\al^{\mu}_a}L_{ij}(z)dz=
\beta^i_a\de_{j\mu} -\de_{i1}\beta_a^{\mu}\al_a^j.
\end{equation}

For $b\neq a$\,, $\al_b$ is a left null vector for the matrix
$||\pa_{\al^{\mu}_a}L_{ij}(\ga_b)||$
$$
\sum_{i=1}^n \al^i_b
\pa_{\al^{\mu}_a}L_{ij}(\ga_b)=0,  \qquad b\neq a,
$$
and the regular part the matrix
$||\pa_{\al^{\mu}_a}L_{ij}(z)||$ at point $\ga_a$ satisfies
the following linear inhomogeneous equation (for the definition of the
matrix $||L^{a,0}_{ij}||$ see equation (\ref{LAX}))

$$
\sum_{i=1}^n \al^i_a \pa_{\al^{\mu}_a}L_{ij}(z) |_{{regular~part~at~}z=\ga_a}=
(k_a\de_{\mu j}-L^{a,0}_{\mu j}).
$$

Finally, the differential $\pa_{\ga_a}L_{ij}(z)dz$ is
holomorphic everywhere on $\Si$ except the point $\ga_a$,
where it has a pole of the second order
and on a neighborhood of the point it
behaves like

\begin{equation}
\label{pagaL}
\pa_{\ga_a}L_{ij}(z)dz=\frac{\beta^i_a \al_a^j dz}{(z-\ga_a)^2} +
{\rm regular~terms}.
\end{equation}

For $b\neq a$\,, $\al_b$ is a left null vector for the matrix
$||\pa_{\ga_a}L_{ij}(\ga_b)||$

$$
\sum_{i=1}^n \al^i_b \pa_{\ga_a}L_{ij}(\ga_b)=0, \qquad b\neq a,
$$
and, in addition, the regular part of the matrix
$||\pa_{\ga_a}L_{ij}(z)||$ at the point $\ga_a$ satisfies
the following linear inhomogeneous
equation (for the definition of the
matrix $||L^{a,1}_{ij}||$ see equation (\ref{LAX}))

$$
\sum_{i=1}^n \al^i_a \pa_{\ga_a}L_{ij}(z) |_{{regular~part~at~}z=\ga_a}=
-\sum_{i=1}^n \al^i_a L^{a,1}_{ij}.
$$

All the properties of the derivatives  $\pa_{k_a}L_{ij}(z)dz$,
$\pa_{\ga_a}L_{ij}(z)dz$,  $\pa_{\al^{\mu}_a}L_{ij}(z)dz$ and
$\pa_{\beta^{\mu}_a}L_{ij}(z)dz$ can be easily derived from
the definition of the Krichever Lax differential (\ref{LAX}) and
the uniqueness of the derivatives as meromorphic
differentials on $\Si$ follows directly from the
Krichever lemma.

\subsection{The sketch of the proof.}
Let us rewrite the Yang-Baxter relation (\ref{LrL}) in
the following form

\begin{equation}
\label{YB}
D_{ijkl}(z,w)dz\otimes dw=R_{ijkl}(z,w)dz\otimes dw,
\end{equation}
where
$$
D_{ijkl}(z,w)=\{L_{ij}(z),L_{kl}(w)\}=
$$
$$
=\sum_{a=1}^{ng}\left( \pa_{\ga_a}L_{ij}(z) \pa_{k_a}L_{kl}(w)-
 \pa_{k_a}L_{ij}(z) \pa_{\ga_a}L_{kl}(w) \right)+
$$
$$
+\sum_{a=1}^{ng}\sum_{\mu=2}^n \left(\pa_{\al^{\mu}_a} L_{ij}(z) \pa_{\beta^{\mu}_a} L_{kl}(w)-
\pa_{\beta^{\mu}_a} L_{ij}(z) \pa_{\al^{\mu}_a} L_{kl}(w)\right),
$$
and
$$
R_{ijkl}(z,w)=\sum_{m=1}^n \de_{il}r_{mk}(z,w) L_{mj}(z)-
L_{il}(z)r_{jk}(z,w)-
$$
$$
-\sum_{m=1}^n \de_{kj} r_{mi}(w,z) L_{ml}(w)+
L_{kj}(w)r_{li}(w,z).
$$

Using the properties of the differentials $\pa_{k_a}L_{ij}(z)dz$,
$\pa_{\ga_a}L_{ij}(z)dz$,  $\pa_{\al^{\mu}_a}L_{ij}(z)dz$ and
$\pa_{\beta^{\mu}_a}L_{ij}(z)dz$ we derive a relatively long
list of properties for the form $D_{ijkl}(z,w)dz\otimes dw$:

\begin{enumerate}

\item The poles of the form $D_{ijkl}(z,w)dz\otimes dw$
are located at the points $\ga_a$ and
$P$ so that the pole at the point $P$ is simple and
the poles at the points $\ga_a$ are of the second order.

\item If $w$ coincides neither
with the point $P$ nor with any of the points $\ga_b$
the singular part of the component
$D_{ijkl}(z,w)$ at the point $z=\ga_a$ looks like

\begin{equation}
\label{Dsing}
D_{ijkl}(z,w)=\frac{D^{a,2}_{ijkl}(w)}{(z-\ga_a)^2}+
\frac{D^{a,1}_{ijkl}(w)}{z-\ga_a}+
{\rm regular~terms},
\end{equation}
where $D^{a,2}_{ijkl}(w)$ is a component of the
holomorphic differential $\al^l_a\beta^i_a\al^j_a u_{ak}(w)dw$
and $D^{a,1}_{ijkl}(w)dw$ is a differential of the
third kind defined by the following properties
\begin{itemize}
\item  $D^{a,1}_{ijkl}(w)dw$ has poles only
at the points $\ga_a$ and $P$ with the
residue at the point $\ga_a$ being

\begin{equation}
\label{Dsing1}
Res_{w=\ga_a} D^{a,1}_{ijkl}(w)dw=
\de_{kj}\beta^i_a\al^l_a-\de_{il}\beta^k_a\al^j_a.
\end{equation}

\item The values of the components
$D^{a,1}_{ijkl}(w)$
at the points $\ga_b,~b\neq a$ satisfy the following
``null vector'' conditions

\begin{equation}
\label{Dsing1'}
\sum_{k=1}^n \al^k_b D^{a,1}_{ijkl}(w) =0, \qquad b\neq a.
\end{equation}

\item The regular parts of $D^{a,1}_{ijkl}(w)$
at the point $\ga_a$ obey the
following linear inhomogeneous equations

\begin{equation}
\label{Dsing1''}
\sum_{k=1}^n \al^k_a D^{a,1}_{ijkl}(w)|_{{\rm regular~part~at~}w=\ga_a}=
-(k_a\de_{il}-L^{a,0}_{il})\al^j_a.
\end{equation}

\end{itemize}

\item The regular parts of the components $D_{ijkl}(z,w)$ at the points
$\ga_a$ satisfy the linear inhomogeneous equations

\begin{equation}
\label{Dreg}
\sum_{i=1}^n \al^i_a D_{ijkl}(z,w)|_{{\rm regular~part~at~}z=\ga_a}=
D^a_{jkl}(w),
\end{equation}
where $D^a_{jkl}(w)$ are components of a meromorphic tensor-valued differential
defined by the following properties

\begin{itemize}
\item $D^a_{jkl}(w)dw$ is holomorphic everywhere on $\Si$ except the
points $P$ and $\ga_a$ where it has poles of the first and
second order respectively.

\item On a neighborhood of the point $w=\ga_a$ it behaves like

\begin{equation}
\label{Dregg}
D^a_{jkl}(w)dw=-\frac{\beta^k_a \al^l_a \al^j_a dw}{(w-\ga_a)^2}+
\frac{(k_a\de_{kj}-L^{a,0}_{kj})\al^l_a dw}{w-\ga_a}+{\rm regular~terms}.
\end{equation}

\item The values of the components $D^a_{jkl}(w)$ at the points $\ga_b,~b\neq a$
satisfy the following ``null vector'' conditions

\begin{equation}
\label{Dregreg}
\sum_{k=1}^n \al^k_b D^a_{jkl}(\ga_b)=0, \qquad b\neq a.
\end{equation}

\item The regular parts of the components $D^a_{jkl}(w)$ at the point
$\ga_a$ obey the following linear inhomogeneous equations

\begin{equation}
\label{Dregregg}
\sum_{k=1}^n \al^k_a D^a_{jkl}(w)|_{{\rm regular~part~at~}w=\ga_a}=
\sum_{k=1}^n
(\al^j_a \al^k_a L^{a,1}_{kl}- \al^l_a \al^k_a L^{a,1}_{kj}).
\end{equation}

\end{itemize}

\end{enumerate}

An analogous detailed analysis of the components $R_{ijkl}(z,w)$ shows
that $R_{ijkl}(z,w)dz\otimes dw$ satisfies all the
above properties of the form $D_{ijkl}(z,w)dz\otimes dw$.
Due to the Krichever lemma
these properties define a unique form $D_{ijkl}(z,w)dz\otimes dw$
and, thus, the  desired statement is proved. $\Box$

\section{Concluding remarks.}
In conclusion, we point out that
the classical $r$-matrix (\ref{simp})
of the extended Krichever Lax matrix (\ref{LAX}) depends only
on the variables $\ga_a$ and $\al_a$, that is,
on coordinates of the respective configuration space.
Since the differential (\ref{LAX}) is linear in the
variables $k_a$ and $\beta^i_a$\,, the genuine
$r$-matrix (\ref{viacont}) of the Hitchin system
also depends only on the variables $\ga_a$ and $\al_a$.

This forces us to assume that the classical
$r$-matrices satisfy simple analogues of classical dynamical
Yang-Baxter equation \cite{EVclass}, which should express the
consistency of the respective Yang-Baxter relations for the
Krichever Lax matrices (\ref{LAX}) and (\ref{laxic}).

Note also that a formal expression for the
classical $r$-matrix of the extended system can be obtained
by the method developed in the
paper \cite{BDOZ}. Following that method we have
to present the system
on the manifold (\ref{Shirokoe}) with
the Krichever Lax matrix
(\ref{LAX}) via an infinite-dimensional
Hamiltonian reduction on $ng$ copies of the
cotangent bundle to the loop group $GL_n(\bbC)[z,z^{-1}]$.
Although the method allows one to express
the desired $r$-matrix in terms of infinite series
in the Krichever-Novikov type basis \cite{KN}, \cite{OK},
\cite{KN1} it turns out to be very hard to analyze
such answers and to identify the resulting $r$-matrix with
any meromorphic object associated with the product of
curves $\Si\times \Si$.

Finally, we mention that it would be interesting to compare
the Krichever parameterization of Lax and $r$-matrix structures
of Hitchin systems based on Tyurin description to the
analogous approach \cite{Schottky} based on the Schottky uniformisation
of Riemann curves and it would be also intriguing to explain
a role of the obtained $r$-matrices in the context of WZNW models
on Riemann surface \cite{KZB}, \cite{KZB1}.

{\bf Acknowledgements.} I would like to express my sincere thanks to
I.M. Krichever and M.A. Olshanetsky for formulating the problem and for
useful discussions of this topic.  I acknowledge I.M. Krichever for constructive
criticisms concerning the first version of this article and A.M. Levin for
an important technical trick, which drastically simplifies the result of this paper.
I also acknowledge H.W. Braden, A.S. Gorsky, S.V. Oblezin and A.V. Zotov
for useful discussions. I am grateful to M. Ching for criticisms
concerning the English language of this paper.
The work is partially supported by
RFBR grant 00-02-17-956, the Grant
for Support of Scientific Schools 00-15-96557,
and the grant INTAS 00-561.

\section*{Appendix. The proof of the Krichever lemma.}
\begin{lem}[Krichever] Let $\nu_i(z)dz$ be
a meromorphic vector-valued differential on the curve $\Si$. Then, for a
generic set of Tyurin parameters $\ga_a\in \Si$ and $\al_a\in \bbC^n$ and for
an arbitrary set of complex numbers $b_a$ there exists a unique meromorphic
vector-valued differential $v_i(z)dz$ having the same singular parts as
the differential $\nu_i(z)dz$ and obeying the following
conditions\footnote{Note that we choose some local coordinates on
neighborhoods of the points $\ga_a$ and the right hand sides of equations
(\ref{lieq}) and (\ref{lieq1}) depend on this choice}:

\begin{itemize}
\item If $v_i(z)dz$ is holomorphic at the point $\ga_a$ then

\begin{equation}
\sum_{i=1}^n v_i(\ga_a)\al^i_a=b_a,
\label{lieq}
\end{equation}

\item and otherwise,

\begin{equation}
\sum_{i=1}^n v_i(z)\al^i_a|_{\rm regular~part~at~z=\ga_a}=b_a.
\label{lieq1}
\end{equation}

\end{itemize}

\end{lem}

~\\
{\bf Proof.}
The statement of the lemma is equivalent to the fact that
for a generic set of Tyurin parameters $(\ga_a,\, \al_a)$ and
for an arbitrary set of complex numbers
$c_a$ there exists a unique holomorphic vector-valued differential $h_i(z)dz$
satisfying the equations

\begin{equation}
\sum_{i=1}^n h_i(\ga_a)\al^i_a=c_a,
\label{lieqq}
\end{equation}
which are, in turn, equivalent to the following
linear inhomogeneous equations

\begin{equation}
\sum_{i=1}^n\sum_{A=1}^g h^{A}_i
\mu_A(\ga_a)\al^i_a=c_a
\label{lieqqq}
\end{equation}
for the expansion coefficients $h^{A}_i$ of the differential
$h_i(z)dz$ in some basis $\{\mu_A(z)dz,~A=1,\ldots, g\}$ of holomorphic
differentials on the curve $\Si$.

Since the number of coefficients $h^{A}_i$ coincides
with the number of equations (\ref{lieqqq}) the
desired statement is equivalent to the fact that the
following $ng\times ng$-matrix

\begin{equation}
M^{(Ai)}_a =\mu_A(\ga_a)\al^i_a
\label{matrix}
\end{equation}
is non-degenerate.

The proof of this fact turns out to be a simple task
of linear algebra.$\Box$

\end{document}